\documentclass{article}

\usepackage{amsmath}
\usepackage{amssymb}
\usepackage{amsthm}
\usepackage{array}
\usepackage{bbold}
\usepackage[T1]{fontenc}
\usepackage[left=2.5cm, right=2.5cm]{geometry}
\usepackage{hyperref}
\usepackage{mdframed}
\usepackage{tikz}
\usepackage{xcolor}

\newtheorem{corollary}{Corollary}
\newtheorem{example}{Example}
\newtheorem{lemma}{Lemma}
\newtheorem{theorem}{Theorem}

\newcommand{\dbar}{\textit{\dj}}
\newcommand{\defeq}{\stackrel{\text{def}}{=}}
\newcommand{\dotp}{\boldsymbol{.}}

\mdfdefinestyle{proof}{bottomline=false,frametitle=Proof,frametitleaboveskip=0,frametitlebelowskip=5pt,
innerleftmargin=3pt,leftmargin=0pt,linewidth=2pt,rightline=false,subtitleaboveline=true,
subtitleaboveskip=5pt,subtitlebelowskip=0pt,subtitleinneraboveskip=1pt,
subtitleinnerbelowskip=0pt,topline=false}

\title{Rate of Convergence of Polynomial Networks to Gaussian Processes}
\author{Adam Klukowski \\ Huawei Noah's Ark Lab \\ \texttt{ak2028@cantab.ac.uk}}
\date{31 October 2021}

\begin{document}

\maketitle

\begin{abstract}
We examine one-hidden-layer neural networks with random weights. It is well-known that in the limit of infinitely many neurons they simplify to Gaussian processes. For networks with a polynomial activation, we demonstrate that the rate of this convergence in 2-Wasserstein metric is $O(n^{-\frac{1}{2}})$, where $n$ is the number of hidden neurons. We suspect this rate is asymptotically sharp. We improve the known convergence rate for other activations, to power-law in $n$ for ReLU and inverse-square-root up to logarithmic factors for erf. We explore the interplay between spherical harmonics, Stein kernels and optimal transport in the non-isotropic setting.
\end{abstract}

\section{Introduction}

We are concerned with a 1-hidden-layer neural network $\mathcal{P}_n$ of width $n$. This is a random function from the sphere $\sqrt{d} S^{d-1} = \big\{ x \in \mathbb{R}^d : \lVert x \rVert^2 = d \big\}$ to $\mathbb{R}$, defined by
\begin{equation*}
\mathcal{P}_n (x) = \tfrac{1}{\sqrt{n}} \sum_{i=1}^n s_i \phi \left( \tfrac{w_i \dotp x}{\sqrt{d}} \right)
\end{equation*}
where $\phi$ is a fixed function $\mathbb{R} \rightarrow \mathbb{R}$ called activation. The randomness comes from the weights $s_i, w_i$, which are random variables taking values in $\mathbb{R}$ and $\mathbb{R}^d$ respectively. We assume they are independent and identically distributed across $i$.

It is well-known that $\mathcal{P}_n$ converges to a Gaussian process (GP) as $n \rightarrow \infty$. This fact was first noticed in~\cite{Neal_Priors}, and discussed in greater generality in~\cite{Hanin_InfiniteNN}. A number of works aimed to go beyond the limit and understand the phenomena in its neighbourhood. The distribution of preactivations\footnote{Intermediate vectors, or neuron outputs, in deep networks comprising multiple stacked affine maps and coordinatewise activations} was studied perturbatively in~\cite{Yaida_NonGaussian}. The behaviour of single observables -- scalars summarizing the distribution, analogous to the moments of a random variable --  was analysed in~\cite{Dyer_Diagrams} and~\cite{Aitken_Asymptotics}. Here we investigate the rate of convergence from the angle of functional metrics, which ``capture the whole distribution''. This line of work was initiated in~\cite{Eldan_Nonasymptotic}, where it was shown that $\mathcal{P}_n$ is within $O \big( n^{-\frac{1}{6}} \big)$ from a GP in the $\infty$-Wasserstein distance when $\phi$ is a polynomial. We improve their method and obtain asymptotically sharp rates in the 2-Wasserstein metric.

This work is closely tied to Stein's method and optimal transport. Stein kernel for $X$ is a matrix-valued function $\tau$ satisfying $\mathbb{E} [X \dotp f(X)] = \mathbb{E} \langle \tau(X), \nabla f \rangle_{HS}$ for any test function $f$. They were introduced in~\cite{Stein_Expectations}; for an overview see for example~\cite{Mijoule_SteinsMethod} or~\cite{Azmoodeh_Survey}. They are strongly related to quantitative and multi-dimensional forms of CLT~\cite{Courtade_Existence}. It was shown in~\cite{Ledoux_HSI} using the Orstein-Uhlenbeck diffusion semigroup that the Wasserstein distance to a Gaussian can be controlled by the discrepancy between the kernel and a constant matrix. Here we show how to exploit the rich symmetries of spherical harmonics to construct Stein kernels. Also, we explore generalizations of diffusive methods~\cite{Ledoux_HSI},\cite{Otto_Generalization} to relate Wasserstein distance and Stein discrepancy in a non-isotropic\footnote{In this paper by isotropic variables we mean those with covariance matrix equal to identity} setting.

\textbf{Organization:} In the rest of this section we define the notation and then describe the main results. Section~\ref{sec_decomposition} relates random functions to vector-valued random variables. We prove our main result about polynomial networks in section~\ref{sec_polynomial}, and analyse general activations in section~\ref{sec_reluerf}. Spherical harmonics are explained in appendix~\ref{apdx_spherical}, and Stein kernels and optimal transport are described in appendix~\ref{apdx_steinkernels}.

\subsection{Notation}

We denote the inner product and induced norm of vectors as
\begin{equation*}
u \dotp v = \sum_{i} u_i v_i \qquad\qquad \lVert u \rVert^2 = u \dotp u
\end{equation*}
Both will often be accompanied by normalizing factors.

We denote the Hilbert-Schmidt product of matrices as
\begin{equation*}
\langle A, B \rangle_{HS} = \text{Tr} \ A B^\top = \sum_{i, j} A_{i, j} B_{i, j}
\end{equation*}

The metric suitable for comparing random objects is the $2$-Wasserstein distance, defined for random vectors $X, Y$ and random functions $f, g$ as
\begin{equation*}
W_2 (X, Y)^2 = \mathbb{E} \left[ \big\lVert X - Y \big\rVert^2 \right] \qquad\qquad W_2 (f, g)^2 = \mathbb{E} \int_{\sqrt{d} S^{d-1}} \big\lvert f(x) - g(x) \big\rvert^2 dx \\
\end{equation*}
When the two vectors do not share a common probability space, we define $\mathcal{W}_2 (X, Y) = \inf_{(X, Y)} W_2 (X, Y)$, where the infimum is taken over all couplings (joint distributions having $X, Y$ as marginals); the definition for random functions is analogous.

Random function $\mathcal{G}$ is a Gaussian process (GP) if the vector $\big( \mathcal{G}(x)\big)_{x \in X}$ has multivariate normal distribution for any finite set of arguments $X$.

We will make heavy use of spherical harmonics -- a set of functions $Y_{l, m} : \sqrt{d} S^{d-1} \rightarrow \mathbb{R}$ indexed by $l \in \mathbb{N}_0, 1 \leq m \leq \dbar_l$. They are discussed in detail in appendix~\ref{apdx_spherical}. Their key property is orthonormality, meaning
\begin{equation*}
\int_{\sqrt{d} S^{d-1}} Y_{l, m} (x) Y_{l', m'} (x) dx = \delta_{l l'} \delta_{m m'}
\end{equation*}
also, they span the Hilbert space of square-integrable functions $\mathcal{L}^2 \big( \sqrt{d} S^{d-1} \big)$. They give rise to the orthogonal family of Gegenbauer polynomials $P_l$.

\subsection{Overview of results and main ideas}

Our main result (with some technicalities omitted) is

\textbf{Theorem~\ref{thm_sphericalkernel}} (simplified)\textbf{.}\textit{
Assume that $s_i$ satisfies $\mathbb{E} [s^2] = 1$, $w_i$ are uniformly distributed on the sphere $w_i \sim \sqrt{d} US^{d-1}$, and the activation $\phi$ is a polynomial. Then there exists a Gaussian process $\mathcal{G}$ such that
\begin{equation*}
\mathcal{W}_2( \mathcal{P}_n, \mathcal{G}) \leq \tfrac{C}{\sqrt{n}}
\end{equation*}
where $C^2 = O \left( (d + \deg \phi)^d \cdot \mathbb{E} \big[ s^4 \big] \cdot \mathbb{E} \Big[ \phi' \big( \mathcal{N} (0, 1) \big)^2 \Big] \right)$.}

The precise statement is theorem~\ref{thm_sphericalkernel} in section~\ref{sec_polynomial}. By approximating the activation with polynomials (section~\ref{sec_reluerf}) we obtain

\textbf{Theorem~\ref{thm_reluerf}} (simplified)\textbf{.}\textit{
Assume $\mathbb{P} (s_i = 1) = \mathbb{P} (s_i = -1) = \tfrac{1}{2}$, $w_i \sim \sqrt{d} U S^{d-1}$. For the rectified linear unit\footnote{$\text{ReLU} (x) = \max \{ 0, x \}$} activation $\phi = \textnormal{ReLU}$ we have
\begin{equation*}
\mathcal{W}_2 \left( \mathcal{P}_n, \mathcal{G} \right) \leq 7 n^{-\frac{3}{2(2d-1)}}
\end{equation*}
while for the error function\footnote{$\text{erf} (x) = \tfrac{2}{\sqrt{\pi}} \int_0^x e^{-y^2} dy$} $\phi = \textnormal{erf}$ we have
\begin{equation*}
\mathcal{W}_2 \left( \mathcal{P}_n, \mathcal{G} \right) \leq C^d (\log n)^\frac{d-2}{2} \cdot n^{-\frac{1}{2}}
\end{equation*}}

The first main idea is to find an embedding $E : \sqrt{d} S^{d-1} \hookrightarrow V$ satisfying $\phi \left( \tfrac{w \dotp x}{\sqrt{d}} \right) = E(w) \dotp E(x)$, and use it to express the neural network as an inner product in the augmented space $V$ as
\begin{equation}\label{eq_nnasproduct}
\mathcal{P}_n (x) = \underbrace{\left( \tfrac{1}{\sqrt{n}} \sum_{i=1}^n s_i E (w_i) \right)}_{\tilde{w}} \dotp E(x)
\end{equation}
This description separates ``random'' from ``function'' -- the first bracket $\tilde{w}$ does not depend on the argument $x$, while the second is a deterministic function of $x$.

The next step is approximating the first factor $\tilde{w}$ by a multivariate normal, using a variant of quantitative CLT. This can then be translated this into an approximation of the network $\mathcal{P}_n$ by a Gaussian process.

In~\cite{Eldan_Nonasymptotic} an embedding into $V = \big( \mathbb{R}^d \big)^{\otimes 0} \oplus\dots\oplus \big( \mathbb{R}^d \big)^{\otimes \deg \phi}$ was obtained by expanding all monomials $(w \dotp x)^k$ of $\phi (w \dotp x)$. Here we use an expansion in the basis of spherical harmonics. This approach gives a simple covariance structure of the random vector $\tilde{w}$ -- its matrix is diagonal with explicit eigenvalues. In fact, this expansion enables us to isometrically translate the problem into a question about countably-dimensional random vectors. Then we employ the machinery of Stein kernels, and construct one by leveraging the geometry of spherical harmonics.

\section{Harmonic decomposition}\label{sec_decomposition}

Here we exhibit a Gaussianity-preserving linear Wasserstein-isometry between random functions and random vectors. Any random function $f$ on the sphere can be expanded in the basis of spherical harmonics, giving a $\mathbb{R}^{\dbar_0 + \dbar_1 + \dots}$-valued random variable
\begin{equation*}
X_{l, m} = \int_{\sqrt{d} S^{d-1}} Y_{l, m}(x) f(x) dx
\end{equation*}
Conversely, every $\mathbb{R}^{\dbar_0 + \dbar_1 + \dots}$-valued random variable $X_{l, m}$ naturally defines a random function via
\begin{equation*}
f(x) = \sum_{l=0}^\infty \sum_{m=1}^{\dbar_l} X_{l, m} Y_{l, m}(x)
\end{equation*}
Spherical harmonics form a complete orthonormal basis, so these transformations are mutually inverse. They define correspondences
\begin{align}
\text{random functions on } \sqrt{d} S^{d-1} \qquad&\longleftrightarrow\qquad \mathbb{R}^{\dbar_0 + \dbar_1 + \dots} \text{-valued random variables} \label{corr_l2norm}\\
\mathbb{E} \int_{\sqrt{d} S^{d-1}} \big\lvert f_1(x) - f_2(x) \big\rvert^2 dx \qquad&\longleftrightarrow\qquad \mathbb{E} \lVert X^{(1)} - X^{(2)} \rVert^2 \label{corr_isometry}\\
\text{Gaussian processeses} \qquad&\longleftrightarrow\qquad \text{multivariate normal variables} \label{corr_gaussianity}\\
\text{NNs } \mathcal{P}_n = \tfrac{1}{\sqrt{n}} \sum_{i=1}^n s_i \phi \left( \tfrac{w_i \dotp x}{\sqrt{d}} \right) \qquad&\longleftrightarrow\qquad \tfrac{1}{\sqrt{n}} \sum_{i=1}^n \tfrac{\hat{\phi}_l}{\sqrt{\dbar_l}} s_i Y_{l, m} (w_i) \label{corr_neural}
\end{align}

Line \ref{corr_isometry} states that we are dealing with an isometry with respect to appropriate 2-Wasserstein metrics. This is a consequence of the orthonormality of harmonics
\begin{align*}
\mathbb{E}& \int_{\sqrt{d} S^{d-1}} \big\lvert f_1(x) - f_2(x) \big\rvert^2 dx = \mathbb{E} \int \left[ \sum_{l, m} \left( X^{(1)}_{l, m} - X^{(2)}_{l, m} \right) Y_{l, m}(x) \right]^2 dx =\\
&= \mathbb{E} \sum_{l, m, l', m'} \left( X^{(1)}_{l, m} - X^{(2)}_{l, m} \right) \left( X^{(1)}_{l', m'} - X^{(2)}_{l', m'} \right) \int Y_{l, m}(x) Y_{l', m'}(x) dx =\\
&= \mathbb{E} \sum_{l, m} \left( X^{(1)}_{l, m} - X^{(2)}_{l, m} \right)^2 = \mathbb{E} \big\lVert X^{(1)} - X^{(2)} \big\rVert^2
\end{align*}

Preservation of Gaussianity~\ref{corr_gaussianity} holds because the maps are linear.

In equation \ref{corr_neural}, the coefficients $\hat{\phi}_l$ come from the expansion $\phi = \sum_{l=0}^\infty \hat{\phi}_l P_l$ of the activation function $\phi$ into Gegenbauer polynomials $P_l$. Equation~\ref{eqgegenbauerpolydefinition} from appendix~\ref{apdx_harmonics} states that their value at a dot product is expressible in terms of spherical harmonics as
\begin{equation*}
P_l \left( \tfrac{w \dotp x}{\sqrt{d}} \right) = \tfrac{1}{\sqrt{\dbar_l}} \sum_{m=1}^{\dbar_l} Y_{l, m}(w) Y_{l, m}(x)
\end{equation*}
this allows us to interpret the network as an Euclidean inner product in an enlarged space
\begin{equation}\label{eq_nnasharmonicproduct}
\mathcal{P}_n (x) = \sum_{l = 0}^\infty \sum_{m=1}^{\dbar_l} \sum_{i=1}^n \tfrac{1}{\sqrt{n}} \cdot \tfrac{\hat{\phi}_l}{\sqrt{\dbar_l}} \cdot s_i Y_{l, m}(w_i) Y_{l, m}(x) = \left( \tfrac{1}{\sqrt{n}} \sum_{i=1}^n \tfrac{\hat{\phi}_l}{\sqrt{\dbar_l}} \cdot s_i Y_{l, m}(w_i) \right)_{l, m} \dotp \Big( Y_{l, m}(x) \Big)_{l, m}
\end{equation}

\section{Polynomial networks}\label{sec_polynomial}

\begin{theorem}\label{thm_sphericalkernel}
Assume that the weights $s_i$ obey $\mathbb{E} [s^2] = 1$, the weights $w_i$ are distributed uniformly on the sphere $w_i \sim \sqrt{d} U S^{d-1}$, and the activation $\phi$ is a polynomial of degree $k$ satisfying $\mathbb{E} \big[ \phi (x_1) \big| x \sim \sqrt{d} U S^{d-1} \big] = 0$. Then, for each $n$, there exists a Gaussian process $\mathcal{G}$ such that
\begin{equation*}
\mathcal{W}_2( \mathcal{P}_n, \mathcal{G}) \leq \tfrac{C}{\sqrt{n}} \qquad\text{where}\qquad C^2 = \tfrac{6 d (d + k)^{d - 2}}{(d - 1)!} \cdot \mathbb{E} \big[ s^4 \big] \mathbb{E} \Big[ \phi' \big( \mathcal{N} (0, 1) \big)^2 \Big] + \textnormal{var} \big[ s^2 \big] \mathbb{E} \Big[ \phi (x_1)^2 \Big\vert x \sim \sqrt{d} U S^{d - 1} \Big]
\end{equation*}
\end{theorem}

\textbf{Idea of proof:} Note it is enough to bound the distance between the random bracket from equation~\ref{eq_nnasharmonicproduct} and a Gaussian. We achieve this by exhibiting a Stein kernel for the random variable $\tfrac{\hat{\phi}_l}{\sqrt{\dbar_l}} Y_{l, m} (w) \ \Big\vert \ w \sim \sqrt{d} U S^{d - 1}$.\\
\begin{minipage}{0.6\textwidth}
Our construction for $l=1$ is illustrated on the right. Recall that $Y_{1, i} (w) = w_i$. For each $i$, we pair up the points $w^+, w^-$ that differ only by the sign of the $i$-th coordinate
\begin{equation*}
\mathbb{E} [w_i f_i(w)] = \tfrac{1}{2} \mathbb{E} \Big[ \lvert w_i \rvert \big( f_i (w^+) - f_i (w^-) \big) \Big]
\end{equation*}
We join them with the shortest curve $\gamma$, and apply the fundamental theorem of calculus to the difference
\begin{equation*}
f_i (w^+) - f_i (w^-) = \int_\gamma \nabla f_i \dotp d \gamma = \mathop{\mathbb{E}}_{w \in \gamma} \Big[ \dot{\gamma}(w) \dotp \nabla f_i (w) \Big]
\end{equation*}
Averaging over the sphere gives an equation of the form
\begin{equation*}
\mathbb{E} [w_i f_i(w)] = \mathbb{E} \big[ (\text{some vector field}) \dotp \nabla f_i \big]
\end{equation*}
Which is precisely the form of a Stein kernel.
\end{minipage}
\hfill
\begin{minipage}{0.37\textwidth}
\centering
\begin{tikzpicture}[thick, scale=0.5]
\draw (0, 0) circle (5);
\draw (-3.28, 4.04) -- (-3, 4) -- (-3.04, 3.72);
\draw (-1.648,  4.936) -- (-1.4, 4.8) -- (-1.536,  4.552);
\draw (-0.2, 5.2) -- (0, 5) -- (-0.2, 4.8);
\draw (1.264, 5.048) -- (1.4, 4.8) -- (1.152, 4.664);
\draw (2.96, 4.28) -- (3, 4) -- (2.72, 3.96);
\filldraw (-4, 3) circle (5pt);
\filldraw (4, 3) circle (5pt);
\node at (-2.5, 1.5) {$w^- \! = \! \begin{pmatrix} \vdots \\ -w_i \\ \vdots \end{pmatrix}$};
\node at (2.5, 1.5) {$\begin{pmatrix} \vdots \\ w_i \\ \vdots \end{pmatrix} \! = \! w^+$};
\node at (0, 4) {$\gamma$};
\end{tikzpicture}
\end{minipage}
It is not immediately clear how to generalize this construction beyond $l=1$. However, it turns out that the vector field we obtain is precisely the gradient of $Y_{1, i}$ tangent to the sphere $\sqrt{d} S^{d-1}$. This interpretation makes sense for any $l, m$. Thus, what we actually do is the calculation of average derivative of test functions in the direction of $\nabla Y_{l, m}$. It turns out that every spherical harmonic except $Y_{l, m}$ is annihilated.

Once we construct the kernel, we compute its Stein discrepancy $S$ using identities from appendix \ref{apdx_rotationsandoperators}. We finish by invoking lemmas from appendix~\ref{apdx_steinkernels} to extract a bound on the Wasserstein distance from the discrepancy.

\begin{mdframed}[style=proof]

\mdfsubtitle{NOTATION AND QUOTED RESULTS}

For $1 \leq l \leq k$ denote
\begin{equation*}
\tilde{Y}_{l, m} (w) = \tfrac{\hat{\phi}_l}{\sqrt{\dbar_l}} Y_{l, m} (w)
\end{equation*}
This is an embedding $\tilde{Y} : \sqrt{d} S^{d - 1} \hookrightarrow \mathbb{R}^{\dbar_1 + \dots + \dbar_k}$, whose left inverse is the projection onto the first $\dbar_1 = d$ coordinates. We will be concerned with the random variable $\tilde{y} = \tilde{Y} (w) \big\vert w \sim \sqrt{d} U S^{d - 1}$.

We will need the rotation matrices $R^\alpha_{ab}$ that act on the basis vectors $e_i$ as
\begin{equation*}
R^\alpha_{ab} e_a = \cos\alpha \ e_a - \sin\alpha \ e_b \qquad\qquad R^\alpha_{ab} e_b = \sin\alpha \ e_a + \cos\alpha \ a_b \qquad\qquad R^\alpha_{ab} e_c = e_c \quad\text{when}\quad c \notin \{ a, b \}
\end{equation*}
and the operators
\begin{equation*}
\partial_r \defeq w \dotp \nabla = \sum_{i=1}^d w_i \partial_i \qquad\qquad L_{a b} \defeq w_a \partial_b - w_b \partial_a \qquad\qquad L^2 \defeq \sum_{a<b} L_{a b}^2
\end{equation*}

Finally we recall the following identities from appendix~\ref{apdx_spherical}
\begin{align}
\big( L_{a b} f \big) (x) =& - \partial_\alpha f( R_{a b}^\alpha x) \Big|_{\alpha=0} \tag{equation~\ref{eq_generators}}\\
r^2 \nabla^2 =& L^2 + \partial_r (\partial_r + d - 2) \tag{equation~\ref{eq_laplacians}}\\
0 =& (d - t^2) P''_l (t) - (d - 1) t P'_l (t) + l (l + d - 2) P_l (t) \tag{equation~\ref{eq_gegenbauerode}}
\end{align}

\mdfsubtitle{CONSTRUCTION OF THE KERNEL}

We want to build a Stein kernel for $\tilde{y}$. Consider a test function $f: \mathbb{R}^{\dbar_1 + \dots + \dbar_k} \rightarrow \mathbb{R}$. We would like to understand
\begin{equation*}
\mathbb{E} \left[ \tilde{y}_{l, m} f(\tilde{y}) \right] = \tfrac{\hat{\phi}_l}{\sqrt{\dbar_l}} \mathbb{E} \left[ Y_{l, m}(w) f(\tilde{Y}(w)) \right]
\end{equation*}
The second expectation is simply the coefficient standing next to $Y_{l, m}$ in the harmonic expansion of $f \circ \tilde{Y}$. We will temporarily move from $\mathbb{R}^{\dbar_1 + \dots + \dbar_k}$ with the test function $f$ to $\mathbb{R}^d \supseteq \sqrt{d} S^{d - 1}$ with the test function $f \circ \tilde{Y}$. As promised, consider the tangent gradient
\begin{equation}\label{eq_pf_tangentgradient}
\nabla Y_{l, m} - \tfrac{l}{d} Y_{l, m} w \qquad \text{and the corresponding operator} \qquad \left( \nabla Y_{l, m} - \tfrac{l}{d} Y_{l, m} w \right) \dotp \nabla
\end{equation}
Viewed in $\mathbb{R}^d$, $Y_{l, m}$ is a homogeneous polynomial of degree $l$, so $w \dotp \nabla Y_{l, m} = l Y_{l, m}$. Hence the operator annihilates $r^2 = \lVert w \rVert^2$, so this vector field is tangent to the sphere $\sqrt{d} S^{d-1}$.

Let us look at how does this vector field affect harmonic expansions. Remembering equation~\ref{eq_laplacians} and $\mathbb{E} [L^2 g] = 0$, we can rewrite the action of~\ref{eq_pf_tangentgradient} on spherical harmonics as
\begin{align*}
\mathbb{E} \Big[ & \left( \nabla Y_{l, m} - \tfrac{l}{d} Y_{l, m} w \right) \dotp \nabla Y_{l', m'} \Big\vert w \sim \sqrt{d} US^{d-1} \Big] = \\
&= \mathbb{E} \Big[ \tfrac{1}{2d} r^2 \nabla^2 \big( Y_{l, m} Y_{l', m'} \big) - \tfrac{l}{d} Y_{l, m} \partial_r Y_{l', m'} \Big] = \\
&= \mathbb{E} \Big[ \tfrac{1}{2d} \big( L^2 + \partial_r (\partial_r + d - 2) \big) \big( Y_{l, m} Y_{l', m'} \big) - \tfrac{l l'}{2 d} Y_{l, m} Y_{l', m'} \Big] =\\
&= \left( \tfrac{(l + l') (l + l' + d - 2)}{2d} - \tfrac{l l'}{d} \right) \mathbb{E} \left[ Y_{l, m} Y_{l', m'} \right] =\\
&= \tfrac{l^2 + l'^2 + (l + l') (d - 2)}{2 d} \delta_{l l'} \delta_{m m'}
\end{align*}

In expectation, this vector field annihilates every spherical harmonic other than $Y_{l, m}$ itself, which is sent to $\tfrac{l (l + d - 2)}{d}$. Therefore we can filter the coefficients of $f \circ \tilde{Y}$ using the identity
\begin{equation*}
\mathbb{E} \left[ Y_{l, m}(w) (f \circ \tilde{Y}) (w) \right] = \tfrac{d}{l (l + d - 2)} \mathbb{E} \left[ \left( \nabla Y_{l, m} - \tfrac{l}{d} Y_{l, m} w \right) \dotp \nabla (f \circ \tilde{Y}) \right]
\end{equation*}

Now we need to return from the sphere $\sqrt{d} US^{d-1}$ and go back to $\mathbb{R}^{\dbar_1 + \dots + \dbar_k}$. We do it using chain rule
\begin{equation*}
\left( \nabla Y_{l, m} - \tfrac{l}{d} Y_{l, m} w \right) \dotp \nabla (f \circ \tilde{Y}) = \sum_{l', m'} \tfrac{\hat{\phi}_l}{\sqrt{\dbar_l}} \left( \nabla Y_{l, m} - \tfrac{l}{d} Y_{l, m} w \right) \dotp \nabla Y_{l', m'} \partial_{l', m'} f
\end{equation*}

Therefore
\begin{align*}
\mathbb{E} \left[ \tilde{y}_{l, m} f(\tilde{y}) \right] =& \mathbb{E} \left[ \sum_{l', m'} \tau_{l, m; l', m'} \partial_{l', m'} f (\tilde{y}) \right] \\
\text{where} \qquad\qquad \tau_{l, m; l', m'} \big( \tilde{Y}(w) \big) =& \tfrac{\hat{\phi}_l \hat{\phi}_{l'}}{\sqrt{\dbar_l \dbar_{l'}}} \tfrac{d}{l (l + d - 2)} \left( \nabla Y_{l, m} - \tfrac{l}{d} Y_{l, m} w \right) \dotp \nabla Y_{l', m'}
\end{align*}
This means that $\tau$ is a Stein kernel for $\tilde{y}$.

\mdfsubtitle{HILBERT-SCHMIDT NORM OF BLOCKS}

Let us rewrite the kernel as
\begin{equation*}
\tau_{l, m; l', m'} = \tfrac{\hat{\phi}_l \hat{\phi}_{l'}}{\sqrt{\dbar_l \dbar_{l'}}} \tfrac{1}{2 l (l + d - 2)} \Big[ r^2 \nabla^2 \left( Y_{l, m} Y_{l', m'} \right) - 2 l l' Y_{l, m} Y_{l', m'} \Big]
\end{equation*}
expanding the Laplacian according to equation~\ref{eq_laplacians} we get
\begin{align*}
r^2 \nabla^2 & \left( Y_{l, m} Y_{l', m'} \right) - l l' Y_{l, m} Y_{l', m'} =\\
=& -(l (l + d - 2) + l' (l' + d - 2)) Y_{l, m} Y_{l', m'} + (l + l') (l + l' + d - 2) Y_{l, m} Y_{l', m'} +\\
&+ \sum_{a, b} L_{a b} Y_{l, m} \cdot L_{a b} Y_{l', m'} - 2 l l' Y_{l, m} Y_{l', m'} =\\
=& \sum_{a, b} L_{a b} Y_{l, m} \cdot L_{a b} Y_{l', m'}
\end{align*}
which means
\begin{equation*}
\tau_{l, m; l', m'} = \tfrac{\hat{\phi}_l \hat{\phi}_{l'}}{\sqrt{\dbar_l \dbar_{l'}}} \tfrac{1}{2 l (l + d - 2)} \sum_{a, b} L_{a b} Y_{l, m} \cdot L_{a b} Y_{l', m'}
\end{equation*}

We will calculate the Hilbert-Schmidt norm of $(l, l')$-block of $\tau$. We have
\begin{align*}
\sum_{m, m'} &\left( \sum_{a, b} L_{a b} Y_{l, m} \cdot L_{a b} Y_{l', m'} \right)^2 =\\
=& \sum_{a, b, c, d, m, m'} L_{a b} Y_{l, m} L_{a b} Y_{l', m'} L_{c d} Y_{l, m} L_{c d} Y_{l', m'} =\\
=& \sum_{a, b, c, d} \left( \sum_m L_{a b} Y_{l, m} L_{c d} Y_{l, m} \right) \left( \sum_{m'} L_{a b} Y_{l', m'} L_{c d} Y_{l', m'} \right)
\end{align*}

We will calculate the sums for fixed $l$. Denote $(S_{a b})_{i j} = \left( \partial_\alpha R^\alpha_{ab} \big|_{\alpha=0} \right)_{i j} = \delta_{a i} \delta_{b j} - \delta_{a j} \delta_{b i}$. Recalling equation~\ref{eq_generators}, we can compute the action of $L$-operators
\begin{align*}
\sum_m & L_{a b} Y_{l, m} L_{c d} Y_{l, m} = \sum_m \partial_\alpha Y_{l, m} \big( R_{a b}^\alpha x \big) \Big\vert_{\alpha = 0} \partial_\beta Y_{l, m} \big( R_{c d}^\beta x \big) \Big\vert_{\beta = 0} =\\
=& \partial_\alpha \partial_\beta \sum_m Y_{l, m} \big( R_{a b}^\alpha x \big) Y_{l, m} \big( R_{c d}^\beta x \big) \Big\vert_{\alpha = \beta = 0} =\\
=& \sqrt{\dbar_l} \partial_\alpha \partial_\beta P_l \Big( \tfrac{1}{\sqrt{d}} x^\top R_{a b}^{-\alpha} R_{c d}^\beta x \Big) \Big\vert_{\alpha = \beta = 0} =\\
=& \tfrac{\sqrt{\dbar_l}}{d} P_l'' (\sqrt{d}) \cdot \partial_\alpha x^\top R_{a b}^{-\alpha} x \Big\vert_{\alpha = 0} \cdot \partial_\beta x^\top R_{c d}^\beta x \Big\vert_{\beta = 0} + \sqrt{\tfrac{\dbar_l}{d}} P_l' (\sqrt{d}) \cdot x^\top \partial_\alpha R_{a b}^{-\alpha} \Big\vert_{\alpha=0} \partial_\beta R_{c d}^\beta \Big\vert_{\beta=0} x =\\
=& - \tfrac{\sqrt{\dbar_l}}{d} P_l'' (\sqrt{d}) \cdot x^\top S_{a b} x \cdot x^\top S_{c d} x - \sqrt{\tfrac{\dbar_l}{d}} P_l' (\sqrt{d}) \cdot x^\top S_{a b} S_{c d} x =\\
=& - \sqrt{\tfrac{\dbar_l}{d}} P_l' (\sqrt{d}) \Big( \delta_{a d} x_b x_c - \delta_{b d} x_a x_c - \delta_{a c} x_b x_d + \delta_{b c} x_a x_d \Big)
\end{align*}

This gives
\begin{align*}
\sum_{m, m'} \tau_{l, m; l', m'}^2 =& \tfrac{\hat{\phi}_l^2 \hat{\phi}_{l'}^2}{\dbar_l \dbar_{l'}} \tfrac{1}{4 l^2 (l + d - 2)^2} \tfrac{\sqrt{\dbar_l \dbar_{l'}}}{d} P_l' (\sqrt{d}) P_{l'}' (\sqrt{d}) \sum_{a, b, c, d} \Big( \delta_{a d} x_b x_c - \delta_{b d} x_a x_c - \delta_{a c} x_b x_d + \delta_{b c} x_a x_d \Big)^2 =\\
=& \tfrac{\hat{\phi}_l^2 \hat{\phi}_{l'}^2}{\sqrt{\dbar_l \dbar_{l'}}} \tfrac{1}{4 d l^2 (l + d - 2)^2} P_l' (\sqrt{d}) P_{l'}' (\sqrt{d}) \cdot 4 (d - 1) r^4
\end{align*}

We substitute $x = \sqrt{d}$ into the differential equation \ref{eq_gegenbauerode} for Gegenbauer polynomials to deduce
\begin{equation*}
P_l'(\sqrt{d}) = \tfrac{l (l + d - 2)}{(d - 1) \sqrt{d}} P_l (\sqrt{d}) = \tfrac{l (l + d - 2)}{(d - 1) \sqrt{d}} \sqrt{\dbar_l}
\end{equation*}

Finally
\begin{equation}\label{eqblocknorm}
\sum_{m, m'} \tau_{l, m; l', m'}^2 = \tfrac{\hat{\phi}_l^2 \hat{\phi}_{l'}^2}{d - 1} \tfrac{l' (l' + d - 2)}{l (l + d - 2)}
\end{equation}

\mdfsubtitle{FINAL BOUND}

Let $\Sigma_{l, m; l', m'} = \delta_{l l'} \delta_{m m'} \tfrac{\hat{\phi}_l^2}{\dbar_l}$. Equation \ref{eqblocknorm} gives
\begin{equation}\label{eq_pf_discrepancy1}
S \Big( \tilde{y}, \Sigma \Big)^2 \leq \big\lVert \Sigma^{-\frac{1}{2}} \tau \big\rVert_{HS}^2 = \sum_{l, l'} \tfrac{\dbar_l}{\hat{\phi}_l^2} \cdot \tfrac{\hat{\phi}_l^2 \hat{\phi}_{l'}^2}{d - 1} \tfrac{l' (l' + d - 2)}{l (l + d - 2)} = \tfrac{1}{d - 1} \left( \sum_{l = 1}^k \tfrac{\dbar_l}{l (l + d - 2)} \right) \left( \sum_{l = 1}^k \hat{\phi}_l^2 l (l + d - 2) \right)
\end{equation}
Substituting $\dbar_l = \binom{d + l - 1}{d - 1} - \binom{d + l - 3}{d - 1}$, for $d \geq 4$ we can bound the first term by
\begin{equation*}
\sum_{l = 1}^k \tfrac{\dbar_l}{l (l + d - 2)} = \sum_{l = 1}^k \tfrac{d^2 + 2dl - 3d - 3l + 2}{l + d - 2} \cdot \tfrac{(d + l - 3) \dots (l + 1)}{(d - 1)!} \leq k \cdot 2 d \cdot \tfrac{(d + k)^{d - 3}}{(d - 1)!} \leq \tfrac{2 d (d + k)^{d - 2}}{(d - 1)!}
\end{equation*}
We can check by hand that this also holds for $d = 2, 3$. The second term of~\ref{eq_pf_discrepancy1} can be simplified by recalling the orthonormality of $P_l$ with respect to the density of single coordinate (equation~\ref{eq_gegenbauerorthonormality} from appendix~\ref{apdx_harmonics})
\begin{equation*}
\int_{-\sqrt{d}}^{\sqrt{d}} P_l (t) P_{l'} (t) \xi (t) dt = \delta_{l l'} \qquad\qquad \text{where} \qquad\qquad \xi (t) = \tfrac{\Gamma \left( \frac{d}{2} \right)}{\Gamma \left( \frac{d-1}{2} \right) \sqrt{\pi d}} \left( 1 - \tfrac{t^2}{d} \right)^{\frac{d-3}{2}}
\end{equation*}
which gives
\begin{equation}\label{eq_sturmliouvillelike}
\sum_{l = 0}^\infty \hat{\phi}_l^2 l (l + d - 2) = \int_{-\sqrt{d}}^{\sqrt{d}} \Bigg( \sum_{l = 0}^\infty \hat{\phi}_l P_l \Bigg) \Bigg( \sum_{l = 0}^\infty \hat{\phi}_l l (l + d - 2) P_l \Bigg) \xi (t) dt
\end{equation}
Recalling the ODE~\ref{eq_gegenbauerode} for Gegenbauer polynomials, we note
\begin{equation*}
\sum_{l = 0}^\infty \hat{\phi}_l l (l + d - 2) P_l = - (d - t^2) \phi'' + (d - 1) t \phi' = -d \left( 1 - \tfrac{t^2}{d} \right)^{-\frac{d - 3}{2}} \left( \left( 1 - \tfrac{t^2}{d} \right)^{\frac{d - 1}{2}} \phi' \right)'
\end{equation*}
substituting this relation and the explicit form of $\xi$ yields
\begin{align*}
\sum_{l = 0}^\infty \hat{\phi}_l^2 l (l + d - 2) =& - \tfrac{\Gamma \left( \frac{d}{2} \right)}{\Gamma \left( \frac{d - 1}{2} \right)} \sqrt{\tfrac{d}{\pi}} \int_{-\sqrt{d}}^{\sqrt{d}} \phi(t) \cdot \left( \left( 1 - \tfrac{t^2}{d} \right)^{\frac{d - 1}{2}} \phi' \right)' dt =\\
=& \tfrac{\Gamma \left( \frac{d}{2} \right)}{\Gamma \left( \frac{d - 1}{2} \right)} \sqrt{\tfrac{d}{\pi}} \int_{-\sqrt{d}}^{\sqrt{d}} \phi'(t)^2 \left( 1 - \tfrac{t^2}{d} \right)^{\frac{d - 1}{2}} dt \leq\\
\leq& \sqrt{\tfrac{d - 1}{2}} \sqrt{\tfrac{d}{\pi}} \int_{-\sqrt{d}}^{\sqrt{d}} \phi'(t)^2 \cdot 2 e^{-\frac{t^2}{2}} dt \leq 2 \sqrt{d (d - 1)} \mathbb{E} \Big[ \phi' \big( \mathcal{N}(0, 1) \big)^2 \Big]
\end{align*}
Finally, equation~\ref{eq_pf_discrepancy1} becomes
\begin{equation*}
S (\tilde{y}, \Sigma)^2 \leq \tfrac{6 d (d + k)^{d - 2}}{(d - 1)!} \cdot \mathbb{E} \Big[ \phi' \big( \mathcal{N} (0, 1) \big)^2 \Big]
\end{equation*}

Now we only need to translate the discrepancy into Wasserstein distance. Lemma~\ref{lemscalingkernel} implies
\begin{equation*}
S (s \tilde{y}, \Sigma)^2 \leq \mathbb{E} \big[ s^4 \big] S (\tilde{y}, \Sigma)^2 + \text{var} \big[ s^2 \big] \cdot \lVert \tilde{y} \rVert^2
\end{equation*}
The norm of $\tilde{y}$ is $\sum_{l=1}^k \hat{\phi}_l^2 \leq \sum_{l=1}^\infty \hat{\phi}_l^2 = \mathbb{E} \Big[ \phi (x_1)^2 \Big\vert x \sim \sqrt{d} U S^{d - 1} \Big]$. By corollary~\ref{corl_clt},
\begin{equation*}
\mathcal{W}_2 \left( \tfrac{1}{\sqrt{n}} \sum_{i=1}^n s_i \tilde{y}_i, \mathcal{N} (0, \Sigma) \right) \leq \tfrac{1}{\sqrt{n}} S \left( s \tilde{y}, \Sigma \right)
\end{equation*}
And according to the theory from section~\ref{sec_decomposition}, this translates isometrically to a distance between $\mathcal{P}_n$ and some Gaussian process.

\hfill $\square$
\end{mdframed}

\section{Non-polynomial activations}\label{sec_reluerf}

Here we obtain approximations of networks with ReLU and erf activations. We do this by truncating the expansion of the activation function into (normalized) Hermite polynomials $h_l$. This is a family of polynomials orthonormal with respect to the Gaussian weight $\tfrac{1}{\sqrt{2 \pi}} e^{-\frac{t^2}{2}}$. They can be viewed as the ``limit'' of Gegenbauer polynomials as $d \rightarrow \infty$. We will make use of the generating function for the normalized Hermite polynomials
\begin{equation}\label{eq_hermitegeneratingfn}
e^{tx - \frac{t^2}{2}} = \sum_{l=0}^\infty \frac{h_l (x) t^l}{\sqrt{l!}}
\end{equation}

\begin{theorem}\label{thm_reluerf}
Assume $d \geq 3$, $s_i \sim U \{ -1, 1 \}$, $w_i \sim \sqrt{d} U S^{d - 1}$. With $\phi = \textnormal{ReLU}$, for sufficiently large $n$ there exists a Gaussian process $\mathcal{G}$ satisfying
\begin{equation*}
\mathcal{W}_2 \left( \mathcal{P}_n, \mathcal{G} \right) \leq 7 n^{-\frac{3}{2(2d-1)}}
\end{equation*}
while with $\phi = \textnormal{erf}$
\begin{equation*}
\mathcal{W}_2 \left( \mathcal{P}_n, \mathcal{G} \right) \leq \sqrt{\tfrac{e}{\log \frac{3}{2}}}^d \cdot \tfrac{(\log n)^\frac{d-2}{2}}{\sqrt{n}}
\end{equation*}
\end{theorem}

\begin{mdframed}[style=proof]

\mdfsubtitle{GENERAL ACTIVATIONS}

Let $\phi = \sum_{l=0}^\infty a_l h_l$ be the expansion of $\phi$ in the basis of normalized Hermite polynomials. Denote the truncations as $\overline{\phi} = \sum_{l=0}^k a_l h_l$ and $\overline{\mathcal{P}}_n (x) = \tfrac{1}{\sqrt{n}} \sum_{i=1} s_i \overline{\phi} \left( \tfrac{w_i \dotp x}{\sqrt{d}} \right)$. Then a simple calculation shows
\begin{equation*}
\mathcal{W}_2 \left( \mathcal{P}_n, \overline{\mathcal{P}}_n \right)^2 = \int_{-\sqrt{d}}^{\sqrt{d}} \big( \phi (t) - \overline{\phi} (t) \big)^2 \xi (t) dt
\end{equation*}
Since $\xi (t) \leq \tfrac{5}{\sqrt{2 \pi}} e^{-\frac{t^2}{2}}$, this is at most
\begin{equation*}
5 \int_{-\infty}^\infty \big( \phi (t) - \overline{\phi} (t) \big)^2 \cdot e^{-\frac{t^2}{2}} \tfrac{dt}{\sqrt{2 \pi}} = 5 \sum_{l=k+1}^\infty a_l^2
\end{equation*}
By theorem~\ref{thm_sphericalkernel}, the truncated network $\overline{\mathcal{P}}_n$ can be approximated by some Gaussian process $\mathcal{G}$ as
\begin{equation*}
\mathcal{W}_2 \left( \overline{\mathcal{P}}_n, \mathcal{G} \right) \leq \sqrt{\tfrac{6 d (d + k)^{d - 2}}{n (d - 1)!} \cdot \mathbb{E} \left[ \phi' (x)^2 \big| x \sim \mathcal{N} (0, 1) \right]}
\end{equation*}
Using the triangle inequality, and simplifying $\tfrac{d (d + k)^{d-2}}{(d-1)!} < \tfrac{d^{d-1} k^{d-2}}{(d-1)!} < e^{d-1} k^{d-2}$, we obtain
\begin{equation}\label{eq_truncationapprox}
\mathcal{W}_2 \left( \mathcal{P}_n, \mathcal{G} \right) \leq \sqrt{5 \sum_{l=k+1}^\infty a_l^2} + \sqrt{\tfrac{3 e^d k^{d - 2}}{n} \cdot \mathbb{E} \left[ \phi' (x)^2 \big| x \sim \mathcal{N} (0, 1) \right]}
\end{equation}

\mdfsubtitle{RELU}

Using the equation~\ref{eq_hermitegeneratingfn}, the coefficients of ReLU satisfy
\begin{equation*}
\sum_{l=0}^\infty \tfrac{a_l t^l}{\sqrt{l!}} = \int_0^\infty x e^{-\frac{(x-t)^2}{2}} \tfrac{dx}{\sqrt{2 \pi}} = \tfrac{e^{-\frac{t^2}{2}}}{\sqrt{2 \pi}} + \tfrac{t}{2} + \tfrac{t}{\sqrt{2 \pi}} \int_0^t e^{-\frac{x^2}{2}} dx = \tfrac{1}{\sqrt{2 \pi}} + \tfrac{t}{2} + \tfrac{1}{\sqrt{2 \pi}} \sum_{l=1}^\infty \tfrac{(-1)^{l-1} t^{2 l}}{l! \cdot 2^l \cdot (2 l - 1)}
\end{equation*}
Which means
\begin{equation*}
a_0 = \tfrac{1}{\sqrt{2 \pi}} \qquad\qquad\qquad a_1 = \tfrac{1}{2} \qquad\qquad\qquad a_l = \tfrac{(-1)^{\frac{l}{2} - 1} \sqrt{l!}}{\sqrt{2 \pi} \cdot \left( \frac{l}{2} ! \right) \cdot 2^{\frac{l}{2}} \cdot (l-1)} \cdot \mathbb{1}_{2 | l} \qquad\text{for}\ l > 1
\end{equation*}
By Stirling's formula $a_l^2 \sim \tfrac{1}{\pi \sqrt{2 \pi}} l^{-\frac{5}{2}}$. Therefore for large enough $l$ we have $a_l^2 < \tfrac{1}{7} l^{-\frac{5}{2}}$, and as a consequence $\sum_{l=k+1}^\infty a_l^2 < \tfrac{1}{7} \sum_{l=k+1}^\infty l^{-\frac{5}{2}} < \tfrac{1}{7} \int_k^\infty l^{-\frac{5}{2}} dl = \tfrac{2}{21} k^{-\frac{3}{2}}$. Inequality~\ref{eq_truncationapprox} becomes
\begin{equation*}
\mathcal{W}_2 \left( \mathcal{P}_n, \mathcal{G} \right) \leq k^{-\frac{3}{4}} + \sqrt{\tfrac{2 e^d k^{d-2}}{n}}
\end{equation*}
Picking $\tfrac{1}{3} n^\frac{2}{2d-1} < k < e^{-1} n^\frac{2}{2d-1}$ makes the two terms be of comparable order, and gives
\begin{equation*}
\mathcal{W}_2 \left( \mathcal{P}_n, \mathcal{G} \right) \leq 3 n^{-\frac{3}{2(2d-1)}} + \sqrt{2 e^d e^{-d+2} n^{-\frac{3}{2d-1}}} < 7 n^{-\frac{3}{2(2d-1)}}
\end{equation*}

\mdfsubtitle{ERF}

Again, using the generating function~\ref{eq_hermitegeneratingfn} we find
\begin{equation*}
\tfrac{\partial}{\partial t} \int_{-\infty}^\infty \text{erf} (x) \cdot e^{tx - \frac{t^2}{2}} \cdot e^{-\frac{x^2}{2}} \tfrac{dx}{\sqrt{2 \pi}} = \sum_{l=0}^\infty \tfrac{lt^{l-1}}{\sqrt{l!}} a_l
\end{equation*}
On the other hand
\begin{align*}
\tfrac{\partial}{\partial t} \int_{-\infty}^\infty \text{erf} (x) \cdot e^{-\frac{(x - t)^2}{2}}\tfrac{dx}{\sqrt{2 \pi}} =& \int_{-\infty}^\infty \text{erf} (x) \cdot \Bigg( - \tfrac{\partial}{\partial x} e^{-\frac{(x - t)^2}{2}} \Bigg) \tfrac{dx}{\sqrt{2 \pi}} = \\
=& -\text{erf} (x) e^{-\frac{(x - t)^2}{2}} \Bigg|_{-\infty}^\infty + \int_{-\infty}^\infty \text{erf}' (x) \cdot e^{-\frac{(x - t)^2}{2}} \tfrac{dx}{\sqrt{2 \pi}} = \\
=& \tfrac{\sqrt{2}}{\pi} \int_{-\infty}^\infty e^{-\frac{3 x^2}{2} + xt - \frac{t^2}{2}} dx = \tfrac{2}{\sqrt{3 \pi}} e^{-\frac{t^2}{3}}
\end{align*}
Comparing the coefficients we obtain $a_l = \tfrac{2 (-1)^\frac{l-1}{2} \sqrt{(l-1)!}}{\sqrt{\pi l} \sqrt{3}^l \left( \frac{l-1}{2} ! \right)} \cdot \mathbb{1}_{2 \nmid l}$.

From Stirling's formula $a_l^2 \sim \left( \tfrac{2}{\pi l} \right)^\frac{3}{2} \left( \tfrac{2}{3} \right)^l$, so eventually $a_l^2 < \left( \tfrac{2}{3} \right)^l$ and $\sum_{l=k+1}^\infty a_l^2 < 2 \left( \tfrac{2}{3} \right)^k$. Then equation~\ref{eq_truncationapprox} together with $\lvert \text{erf}' \rvert \leq 1$ give
\begin{equation*}
\mathcal{W}_2 \left( \mathcal{P}_n, \mathcal{G} \right) \leq \sqrt{10} \left( \tfrac{2}{3} \right)^\frac{k}{2} + \sqrt{\tfrac{3 e^d k^{d-2}}{n}}
\end{equation*}
Setting $k \sim \tfrac{\log n}{\log \frac{3}{2}}$ completes the proof.

\hfill $\square$
\end{mdframed}

\section{Discussion}

We have demonstrated that one-hidden-layer neural networks with polynomial activation approach GPs at the rate $O( n^{-\frac{1}{2}} )$ in 2-Wasserstein distance. A natural question to ask is how far can our result be generalized. Can the condition of a polynomial activation be dropped? How about $p$-Wasserstein metrics for $p>2$? The difficulty with the former question seems to originate from $\Sigma^{-1}$ in the definition of $S$ in lemma~\ref{lem_distancefromdiscrepancy}. This factor does not appear in the isotropic case, and for general covariances one may define $S$ in a few different ways and still obtain bounds on $\mathcal{W}$. However, dimensional analysis suggests that $\Sigma^{-1}$ is more than just an artifact of a particular wording of Cauchy-Schwarz inequality: if we scale every variable by $\lambda$, then $\mathcal{W} (X, \mathcal{N})^2$ scales like $\lambda^2$ but plain $\lVert \tau - \Sigma \rVert_{HS}^2$ scales like $\lambda^4$. In our proof, construction of the kernel and relating the distance to discrepancy are largely independent, so we hope that deeper understanding of the relationship between Wasserstein distance and Stein discrepancy will allow to improve our result for little extra effort.

In classical CLT the convergence of a normalized sum to a Gaussian is not faster than $O( n^{-\frac{1}{2}} )$, provided that the variables being averaged have non-zero fourth cumulant. Therefore the bound from theorem~\ref{thm_sphericalkernel} is likely to be asymptotically sharp in $n$.

\appendix

\section{Spherical harmonics}\label{apdx_spherical}

There are at least three equivalent ways to think about spherical harmonics
\begin{itemize}
\item Algebraic: harmonic (i.e. $\nabla^2 Y = 0$) homogeneous polynomials in $d$ variables
\item Representation-theoretic: irreducible representations of $SO(d)$
\item Analytic: basis of the Hilbert space $\mathcal{L}^2 (\sqrt{d} S^{d - 1})$ of functions on the sphere
\end{itemize}

In the discussion below, we will typically start with the algebraic picture, viewing polynomials as elements of the ring $\mathbb{R} [X_1, \dots, X_d]$ and operators $L_{a b}$ as $\mathbb{R}$-linear derivations\footnote{Operator $L$ is a derivation if it satisfies the Leibniz rule $L (fg) = L f \cdot g + f \cdot L g$} over this ring. Then we proceed to the analytic picture -- remind ourselves that polynomials can be treated as functions on the sphere $\sqrt{d} S^{d - 1}$, and think of $L_{a b}$ as infinitesimal generators of rotations; we translate the algebraic results and explore the consequences of acquiring an inner product. Representation-theoretic picture will be present in the background and manifest itself whenever we talk about the symmetries of spherical harmonics.

\subsection{Rotations and operators}\label{apdx_rotationsandoperators}

Special orthogonal group $SO(d)$ acts on points from $\sqrt{d} S^{d-1}$ by $R : x \mapsto Rx$, preserving geometry $x_1 \dotp x_2$. This induces an action on functions $\mathcal{L}^2 (\sqrt{d} S^{d-1})$ by $R : f \mapsto f \circ R^{-1}$, which preserves the inner product
\begin{equation}\label{eq_functionproduct}
\int_{\sqrt{d} S^{d - 1}} f_1(x) f_2(x) dx = \mathbb{E} \big[ f_1(x) f_2(x) \big\vert x \sim \sqrt{d} U S^{d - 1} \big] \qquad\text{with normalization } \int 1 dx = 1
\end{equation}
The group is generated by rotations of the form $R_{a b}^\alpha = \exp (\alpha S_{a b})$ for anti-symmetric matrices $(S_{a b})_{i j} = \delta_{a i} \delta_{b j} - \delta_{a j} \delta_{b i}$. They act on the basis vectors as
\begin{equation*}
R^\alpha_{ab} e_a = \cos\alpha \ e_a - \sin\alpha \ e_b \qquad\qquad R^\alpha_{ab} e_b = \sin\alpha \ e_a + \cos\alpha \ a_b \qquad\qquad R^\alpha_{ab} e_c = e_c \quad\text{when}\quad c \notin \{ a, b \}
\end{equation*}
The infinitesimal generators of such rotations are
\begin{equation}
\partial_\alpha R_{a b}^\alpha f \Big\vert_{\alpha = 0} = L_{a b} f \qquad\text{where}\qquad L_{a b} = X_a \partial_b - X_b \partial_a \label{eq_generators}
\end{equation}
In particular, $\mathbb{E} [L_{a b} f (x) | x \sim \sqrt{d} U S^{d-1}] = \partial_\alpha \mathbb{E} [R_{a b}^\alpha f] \big|_{\alpha=0} = 0$.

Define the Laplace-Beltrami operator as
\begin{equation*}
L^2 \defeq \sum_{a < b} L_{a b}^2 = \tfrac{1}{2} \sum_{a, b} L_{a b}^2
\end{equation*}
it is straightforward to verify
\begin{equation}\label{eq_laplacians}
r^2 \nabla^2 = L^2 + \partial_r (\partial_r + d - 2)
\end{equation}
where
\begin{equation*}
\partial_r \defeq x \dotp \nabla = \sum_{i=1}^d X_i \partial_i
\end{equation*}
Let us note the algebraic properties of these operators. The simplest is $\partial_r$ -- it multiplies a polynomial by its degree. Operators $L_{a b}$ are derivations annihilating $r^2 \defeq X_1^2 + \dots + X_d^2$, so Laplace-Beltrami operator satisfies $L^2 (r^2 f) = r^2 L^2 f$, and by identity \ref{eq_laplacians} is multiplies harmonic homogeneous polynomials of degree $l$ by $-l (l + d - 2)$. Both $\partial_r$ and $L^2$ are invariant under rotations (equivalently, commute with each $L_{a b}$).

Now we describe their basic analytic properties. The operators $L_{a b}$ are tangent\footnote{In the sense $L_{a b} r^2 = 0$, or $L_{a b} = v \dotp \nabla$ where the vector field $v$ is tangent to $\sqrt{d} S^{d-1}$} to the sphere, so they and $L^2$ have well-defined restrictions\footnote{Any derivation $J$ on $\mathbb{R} [X_1, \dots, X_d]$ annihilating $r^2$ gives rise to an operator on $\mathcal{L}^2 (\sqrt{d} S^{d-1})$ as follows. For a function $f$ that is the restriction of a polynomial $F$ to $\sqrt{d} S^{d-1}$ we send $f \mapsto J(F) \big|_{\sqrt{d} S^{d-1}}$. This is well defined, because if $F \big|_{\sqrt{d} S^{d-1}} = F' \big|_{\sqrt{d} S^{d-1}}$ then $F - F' = (r^2 - d) G$ for some $G$, so $J(F - F')$ is the zero function on $\sqrt{d} S^{d-1}$. Conversely, a differential operator $\sum_{i=1}^d P_i \partial_i$ acting on $\mathcal{L}^2 (\sqrt{d} S^{d-1})$ with $P_i \in \mathbb{R} [X_1, \dots, X_d]$ can be naturally reinterpreted as a derivation on $\mathbb{R} [X_1, \dots, X_d]$.} to $\mathcal{L}^2 (\sqrt{d} S^{d - 1})$. Inner product of functions is invariant under rotations and $L_{a b}$ obey Leibniz rule, so they are anti-self-adjoint
\begin{equation*}
\int_{\sqrt{d} S^{d - 1}} L_{a b} f_1 \cdot f_2 dx = - \int_{\sqrt{d} S^{d - 1}} f_1 \cdot L_{a b} f_2 dx
\end{equation*}
As a consequence, $L^2$ is self-adjoint (with respect to~\ref{eq_functionproduct}). Also, $L_{a b}$ annihilate the constant function, so again $\mathbb{E} [L_{a b} f (x) | x \sim \sqrt{d} U S^{d-1} ] = 0$.

Finally let us note a geometric fact about the sphere. It will be usefult later to know the distribution of a single coordinate $x_1$ when we draw $x$ uniformly from the sphere $\sqrt{d} S^{d-1}$. Its density is supported on the interval $[-\sqrt{d}, \sqrt{d}]$ and equals
\begin{equation}\label{eq_singledensity}
\xi (x) = \tfrac{\Gamma \left( \frac{d}{2} \right)}{\Gamma \left( \frac{d-1}{2} \right) \sqrt{\pi d}} \left( 1 - \tfrac{x^2}{d} \right)^\frac{d - 3}{2}
\end{equation}
One way to see this is by noting that $\tfrac{x_1^2}{d}$ for $x \sim \sqrt{d} U S^{d-1}$ has the same distribution as $\tfrac{z_1^2}{\lVert z \rVert^2}$ for $z \sim \mathcal{N} (0, I_d)$, that is $\text{B} \left( \tfrac{1}{2}, \tfrac{d-1}{2} \right)$.

\subsection{Spherical harmonics}\label{apdx_harmonics}

\begin{lemma}\label{lemharmonicdecomposition}
Every homogeneous polynomial $f \in \mathbb{R} [X_1, \dots, X_d]$ can be uniquely written as
\begin{equation*}
f = f_0 + r^2 f_1 + r^4 f_2 + \dots
\end{equation*}
where $f_i$ are homogeneous harmonic polynomials.
\end{lemma}
\begin{mdframed}[style=proof]
We proceed by induction on $l = \deg f$. For $l = 0, 1$ the statement is trivial. For $l \geq 1$, by inductive assumption we may write
\begin{equation*}
\nabla^2 f = g_0 + r^2 g_1 + \dots + r^{2 \left\lfloor \frac{l - 2}{2} \right\rfloor} g_{\left\lfloor \frac{l - 2}{2} \right\rfloor}
\end{equation*}
for harmonic $g_i$ of degree $l - 2 - 2 i$. Now, construct
\begin{equation*}
g = \sum_{i=0}^{\left\lfloor \frac{l - 2}{2} \right\rfloor} \frac{r^{2i} g_i}{2 (i + 1) (2 l - 2 i - 4 + d)}
\end{equation*}
then, either by writing $\nabla^2 (r^2 \bullet) = L^2 + \partial_r^2 + (d + 2) \partial_r + 2d$ and recalling the eigenvalues of $L^2$ and $\partial_r$, or by direct calculation, we can see that
\begin{equation*}
\nabla^2 (r^2 g) = \sum_{i = 0}^{\left\lfloor \frac{l - 2}{2} \right\rfloor} r^{2 i} g_i = \nabla^2 f
\end{equation*}
Therefore, the Laplacian of $f - r^2 g$ is zero, and $f$ can be decomposed as
\begin{equation*}
f = \underbrace{\big( f - r^2 g \big)}_{f_0} + \sum_{i=0}^{\left\lfloor \frac{l - 2}{2} \right\rfloor} r^{2 i + 2} \frac{g_i}{2 (i + 1) (2 l - 2 i - 4 + d)}
\end{equation*}

To see uniqueness, note that each factor is an eigen-element of $L^2$ with a different eigenvalue.

\hfill $\square$
\end{mdframed}

Let us denote the space of degree $l$ homogeneous harmonic polynomials as $H_l$. By the lemma \ref{lemharmonicdecomposition} above, we have
\begin{equation*}
\{ \deg \text{-} l \text{ homog polys} \} = H_l \ \oplus \ r^2 \{ \deg \text{-} (l - 2) \text{ homog polys} \} = H_l \oplus r^2 H_{l - 2} \oplus r^4 H_{l - 4} \oplus \dots
\end{equation*}
this allows to deduce their dimensions
\begin{equation*}
\dbar_l \defeq \dim H_l = \binom{d + l - 1}{d - 1} - \binom{d + l - 3}{d - 1}
\end{equation*}
Harmonicity and homogeneity of given degree are preserved by rotations, so each $H_l$ is closed under $SO(d)$, and each $r^{2 k} H_l$ is a subrepresentation of $SO(d)$ inside $\mathbb{R} [X_1, \dots, X_d]$. Note that by equation \ref{eq_laplacians} we have $H_l = \ker L^2 + l (l + d - 2)$ (in algebraic sense, with $X_i$ considered as abstract symbols).

Now let us think about restrictions of polynomials $\mathbb{R} [X_1, \dots, X_n]$ to $\sqrt{d} S^{d - 1}$. By Stone-Weierstrass theorem, they are dense in $C (\sqrt{d} S^{d - 1}, \mathbb{R})$ with supremum norm. Thus with $\ell_2$-norm we must have
\begin{equation*}
\mathcal{L}^2 \big( \sqrt{d} S^{d-1} \big) = \overline{\bigoplus_{l=0}^\infty H_l}
\end{equation*}
with each $H_l$ closed under $SO(d)$. Also, $H_l = \ker L^2 + l (l + d - 2)$ (in analytic sense, with $H_l$ considered as functions on the sphere and $L^2$ as a second-order differential operator); the operator $L^2$ is self-adjoint, so different $H_l$ are orthogonal.

We take spherical harmonics $Y_{l, 1}, \dots, Y_{l, \dbar_l}$ to be any orthonormal basis of $H_l$. Then
\begin{equation*}
\mathbb{E} \big[ Y_{l, m} (w) Y_{l', m'} (w) \big\vert w \sim \sqrt{d} U S^{d - 1} \big] = \delta_{l l'} \delta_{m m'}
\end{equation*}
Each $H_l$ comes with a representation $\rho$ of $SO(d)$
\begin{equation*}
R Y_{l, m} = Y_{l, m} (R^{-1} x) = \sum_{m' = 1}^{\dbar_l} \rho(R)_{m, m'} Y_{l, m'} (x)
\end{equation*}
Such matrices $\rho(R)$ are also orthogonal, which follows from the invariance of the inner product:
\begin{align*}
\big( \rho(R) \rho(R)^\top \big)_{m, m''} =& \sum_{m'} \rho(R)_{m, m'} \rho(R)_{m'', m'} = \sum_{m', m'''} \rho(R)_{m, m'} \rho(R)_{m'', m'''} \big\langle Y_{l, m'}, Y_{l, m'''} \big\rangle =\\
=& \big\langle R Y_{l, m}, R Y_{l, m''} \big\rangle = \big\langle Y_{l, m}, Y_{l, m''} \big\rangle = \delta_{m, m''}
\end{align*}

Now we look at the relation between spherical harmonics at different points, which will eventually lead to Gegenbauer polynomials. Consider
\begin{equation*}
\check{P}_{l, x} (x') \defeq \tfrac{1}{\sqrt{\dbar_l}} \sum_{m = 1}^{\dbar_l} Y_{l, m} (x) Y_{l, m} (x')
\end{equation*}
By construction $\check{P}_{l, x} \in H_l$. Also, for any rotation $R \in SO(d)$ we have
\begin{equation*}
\check{P}_{l, Rx} (R x') = \tfrac{1}{\sqrt{\dbar_l}} Y_{l, :} (R x)^\top Y_{l, :} (R x') = \tfrac{1}{\sqrt{\dbar_l}} Y_{l, :} (x)^\top \rho(R^{-1})^\top \rho(R^{-1}) Y_{l, :} (x') = \tfrac{1}{\sqrt{\dbar_l}} Y_{l, :} (x)^\top Y_{l, :} (x') = \check{P}_{l, x} (x')
\end{equation*}
Therefore $\check{P}_{l, x} (x')$ depends only on the angle between $x, x'$ and not on their absolute position on the sphere, i.e. $P_{l, x} (x') = P_l \left( \tfrac{x \dotp x'}{\sqrt{d}} \right)$ for some function $P_l$; it must be a polynomial of degree at most $l$. This gives us the key identity
\begin{equation}\label{eqgegenbauerpolydefinition}
P_l \left( \tfrac{x \dotp x'}{\sqrt{d}} \right) = \tfrac{1}{\sqrt{\dbar_l}} \sum_{m = 1}^{\dbar_l} Y_{l, m} (x) Y_{l, m} (x')
\end{equation}
The $P_l$ are called Gegenbauer polynomials\footnote{Different scaling/normalization conventions are used the literature}. They are the unique (up to scaling) functions for which the map $x \mapsto P_l (x_1)$ belongs to $H_l$. Orthogonality of the spaces $H_l$ for different $l$ means that Gegenbauer polynomials are orthogonal with respect to the single-coordinate density $\xi$ from equation~\ref{eq_singledensity}
\begin{align}
P_l (\sqrt{d}) =& \tfrac{1}{\sqrt{\dbar_l}} \sum_{m = 1}^{\dbar_l} \mathbb{E} \left[ Y_{l, m} (x) Y_{l, m} (x) \right] = \sqrt{\dbar_l} \nonumber \\
\begin{split}\label{eq_gegenbauerorthonormality}
\int_{- \sqrt{d}}^{\sqrt{d}} P_l (t) & P_{l'} (t) \cdot \xi (t) dt = \mathbb{E} \left[ P_l (X_i) P_{l'} (X_i) \right] =\\
=& \tfrac{1}{\sqrt{\dbar_l \dbar_{l'}}} \sum_{m, m'} Y_{l, m} (\sqrt{d} e_i) Y_{l', m'} (\sqrt{d} e_i) \mathbb{E} \left[ Y_{l, m} (x) Y_{l', m'} (x) \right] = \tfrac{1}{\sqrt{\dbar_l}} \delta_{l l'} P_l (\sqrt{d}) = \delta_{l l'}
\end{split}
\end{align}
therefore $P_l$ can be computed by Gram-Schmidt orthonormalization of $\{ t^0, t^1, t^2, \dots \}$ with respect to the the density of a single coordinate $\xi$.

Finally we exhibit an ODE for $P_l$. Observe that $r^l P_l \left( \tfrac{X_1 \sqrt{d}}{r} \right)$ is a homogeneous degree-$l$ harmonic polynomial. After tidying up the harmonicity condition we obtain
\begin{equation}\label{eq_gegenbauerode}
0 = \tfrac{1}{r^{l - 2}} \nabla^2 \left( r^l P_l \left( \tfrac{ X_1 \sqrt{d}}{r} \right) \right) = (d - t^2) P''_l (t) - (d - 1) t P'_l (t) + l (l + d - 2) P_l (t)
\end{equation}

\begin{example}
For $d = 2$ this construction is precisely the Fourier analysis. We work over $\sqrt{2} S^1 = \{ (x_1, x_2) : x_1^2 + x_2^2 = 2 \}$, parameterized as $x_1 = \sqrt{2} \cos \theta, x_2 = \sqrt{2} \sin \theta$. Harmonic subspaces are
\begin{align*}
H_0 =& \ \textnormal{span} \big\{ Y_{0, 1} = 1 \big\} \quad&\text{with}\qquad \dbar_0 = 1 \\
H_l =& \ \textnormal{span} \big\{ Y_{l, 1} = \sqrt{2} \cos l \theta, Y_{l, 2} = \sqrt{2} \sin l \theta \big\} \quad&\text{with}\qquad \dbar_l = 2
\end{align*}
spherical harmonics are restrictions of polynomials
\begin{equation*}
Y_{l, 1} = 2^{\frac{1 - l}{2}} \Re (X_1 + i X_2)^l = \tfrac{r^l}{\sqrt{2}^{l - 1}} \cos l \theta \qquad\qquad Y_{l, 2} = 2^{\frac{1 - l}{2}} \Im (X_1 + i X_2)^l = \tfrac{r^l}{\sqrt{2}^{l - 1}} \sin l \theta
\end{equation*}
There is only one rotation generator
\begin{align*}
L_{1 2} =& \partial_\theta = X_1 \partial_2 - X_2 \partial_1\\
L^2 =& \partial^2_\theta = X_1^2 \partial_2^2 + X_2^2 \partial_1^2 - 2 X_1 X_2 \partial_1 \partial_2 - X_1 \partial_1 - X_2 \partial_2
\end{align*}
and the Laplace-Beltrami operator $L^2$ acts on $H_l$ as a multiplication by $-l^2$.

Gegenbauer polynomials are characterized by
\begin{equation*}
P_l \big( \sqrt{2} \cos (\theta - \theta') \big) = \sqrt{2} \cos l \theta \cos l \theta' + \sqrt{2} \sin l \theta \sin l \theta' = \sqrt{2} \cos l (\theta - \theta')
\end{equation*}
i.e. are rescaled Chebyshev polynomials. They are orthonormal with respect to $\xi(t) = \tfrac{dt}{\pi \sqrt{2 - t^2}} = \tfrac{d \theta}{\pi}$.
\end{example}

\section{Stein kernels}\label{apdx_steinkernels}

We say that $\tau$ is a Stein kernel for random variable $X$ if for each $f \in C^\infty_c$ we have
\begin{equation*}
\mathbb{E} \big[ X \dotp f(X) \big] = \mathbb{E} \big\langle \tau(X), \text{Jac} f (X) \big\rangle_{HS}
\end{equation*}
where $(\text{Jac} f)_{a b} = \tfrac{\partial f_a}{\partial X_b}$ is the Jacobian of $f$, and $\langle A, B \rangle_{HS} = \text{Tr} \ A B^\top$ is the Hilbert-Schmidt product.

One can show that a constant matrix $\Sigma$ is a Stein kernel for $X$ if and only if $X \sim \mathcal{N}(0, \Sigma)$ (this statement is known as Stein's lemma). It turns out that the difference between $\tau$ and $\Sigma$ can be used to bound the Wasserstein distance between $X$ and $\mathcal{N} (0, \Sigma)$ (see lemma \ref{lem_distancefromdiscrepancy}). The measure of deviation is called Stein discrepancy, and in the isotropic case it is defined as $S(X, \text{Id}) = \inf_\tau \mathbb{E} \lVert \tau(X) - \text{Id} \rVert_{HS}^2$. We will be working with non-isotropic random variables, and following the formulation of lemma \ref{lem_distancefromdiscrepancy} we generalize the Stein discrepancy as
\begin{equation*}
S (X, \Sigma) \defeq \inf_\tau \sqrt{\mathbb{E} \big\lVert \Sigma^{-\frac{1}{2}} (\tau(X) - \Sigma) \big\rVert_{HS}^2}
\end{equation*}
However, note that other generalizations to non-isotropic case are also possible, and modifying the last part of the proof\footnote{For example by rearranging the equation~\ref{eq_flownorm_raw} before applying Cauchy-Schwarz inequality} of \ref{lem_distancefromdiscrepancy} can give bounds of a different form.

Substituting $f(X) = X_i e_j$ we see that $\mathbb{E} \tau = \mathbb{E} X X^\top = \text{cov} [X]$. Therefore, Stein discrepancy can also be viewed as a measure of variance of $\tau$. Intuitively, as we average independent copies of $X$, we can expect the variance to decrease and $\tau$ to approach its expectation, leading to central limit theorem. This intuition is formalized in corollary \ref{corl_clt}; a stronger result -- that $n S \left( \tfrac{1}{\sqrt{n}} \sum_{i=1}^n X_i \right)$ is non-increasing in $n$ -- was proved in~\cite{Courtade_Existence}.

\subsection{Addition and scaling}

\begin{lemma}\label{lemadditivekernel}
Suppose $\tau_1, \dots, \tau_n$ are Stein kernels for independent $X_1, \dots, X_n$, and write $\bar{X} = \sum_{i=1}^n X_i$. Then
\begin{equation*}
\tau(x) \defeq \mathbb{E} \left[ \sum_i \tau_i(X_i) \Bigg\vert \bar{X} = x \right]
\end{equation*}
is a Stein kernel for $\bar{X}$. If $X_i$ have the same covariance $\Sigma$, then $S \big( \bar{X}, n \Sigma \big)^2 \leq \tfrac{1}{n} \sum_{i = 1}^n S \big( X_i, \Sigma \big)^2$.
\end{lemma}
\begin{mdframed}[style=proof]
\begin{equation*}
\mathbb{E} \left[ \bar{X} \dotp f \big( \bar{X} \big) \right] = \sum_i \mathbb{E} \left[ X_i \dotp f \big( \bar{X} \big) \right] = \sum_i \mathbb{E} \left\langle \tau_i \big( X_i \big), \text{Jac} f \big( \bar{X} \big) \right\rangle_{HS} = \mathop{\mathbb{E}}_{\bar{X}} \left\langle \mathbb{E} \left[ \sum_i \tau_i (X_i) \Bigg\vert \bar{X} \right], \text{Jac} f \big( \bar{X} \big) \right\rangle_{HS}
\end{equation*}
If all covariances are equal then we have
\begin{align*}
S \big( \bar{X}, n \Sigma \big)^2 =& \mathop{\mathbb{E}}_{\bar{X}} \Bigg\lVert \big( n \Sigma \big)^{-\frac{1}{2}} \Bigg( \mathbb{E} \Bigg[ \sum_i \tau_i \big( X_i \big) \Bigg\vert \bar{X} \Bigg] - n \Sigma \Bigg) \Bigg\rVert_{HS}^2 = \tfrac{1}{n} \mathop{\mathbb{E}}_{\bar{X}} \Bigg\lVert \mathbb{E} \Bigg[ \sum_i \Sigma^{-\frac{1}{2}} \Big( \tau_i \big( X_i \big) - \Sigma \Big) \Bigg\vert \bar{X} \Bigg] \Bigg\rVert_{HS}^2 \leq\\
\leq& \tfrac{1}{n} \mathbb{E} \Bigg\lVert \sum_i \Sigma^{-\frac{1}{2}} \Big( \tau_i \big( X_i \big) - \Sigma \Big) \Bigg\rVert_{HS}^2 = \tfrac{1}{n} \sum_i \mathbb{E} \Big\lVert \Sigma^{-\frac{1}{2}} \Big( \tau_i \big( X_i \big) - \Sigma \Big) \Big\rVert_{HS}^2
\end{align*}

\hfill $\square$
\end{mdframed}

\begin{lemma}\label{lemscalingkernel}
If $\tau$ is a Stein kernel for $X$ and $Y = s X$, then $\tau' (y) = \mathbb{E} \big[ s^2 \tau (X) \big\vert Y = y \big]$ is a Stein kernel for $Y$. Its discrepancy is at most
\begin{equation*}
S \big( s X, \mathbb{E} [s^2] \textnormal{cov} [X] \big)^2 \leq \tfrac{\mathbb{E} s^4}{\mathbb{E} s^2} \cdot S (X, \textnormal{cov} [X])^2 + \tfrac{\textnormal{var} [s^2]}{\mathbb{E} s^2} \cdot \mathbb{E} \lVert X \rVert^2
\end{equation*}
\end{lemma}
\begin{mdframed}[style=proof]
\begin{align*}
\mathbb{E} \big[ s X \dotp f (s X) \big] =& \mathop{\mathbb{E}}_X \big[ X \dotp \mathop{\mathbb{E}}_s \big[ s f (s X) \big] \big] = \mathop{\mathbb{E}}_X \Big\langle \tau (X), \text{Jac} \mathop{\mathbb{E}}_s \big[ s f (s X) \big] \Big\rangle_{HS} =\\
=& \mathop{\mathbb{E}}_X \Big\langle \tau(X), \mathop{\mathbb{E}}_s \big[ s^2 \big( \text{Jac} f \big) (s X) \big] \Big\rangle_{HS} = \mathop{\mathbb{E}}_{X, s} \Big\langle s^2 \tau(X), \big( \text{Jac} f \big) (s X) \Big\rangle_{HS} =\\
=& \mathop{\mathbb{E}}_Y \Big\langle \mathbb{E} \big[ s^2 \tau (X) \big\vert Y \big], \text{Jac} f (Y) \Big\rangle_{HS} = \mathbb{E} \Big\langle \tau'(Y), \text{Jac} f (Y) \Big\rangle_{HS}
\end{align*}

Now we will bound its discrepancy. Denote $\mathbb{E} s^2 = \sigma^2, \text{cov} [X] = \Sigma$. Then
\begin{align*}
S \big( s X, \mathbb{E} \big[ s^2 \big] \text{cov} [X] \big)^2 \leq & \mathop{\mathbb{E}}_Y \Big\lVert \mathbb{E} \Big[ \sigma^{-1} \Sigma^{-\frac{1}{2}} \big( s^2 \tau(X) - \sigma^2 \Sigma \big) \Big\vert Y \Big] \Big\rVert_{HS}^2 \leq \mathbb{E} \Big\lVert \sigma^{-1} \Sigma^{-\frac{1}{2}} \big( s^2 \tau(X) - \Sigma \big) \Big\rVert_{HS}^2 =\\
=& \sigma^{-2} \mathbb{E} \big[ s^4 \big] \mathbb{E} \Big\lVert \Sigma^{-\frac{1}{2}} \tau(X) \Big\rVert_{HS}^2 - \sigma^{-2} \big\lVert \Sigma^{\frac{1}{2}} \big\rVert_{HS}^2 = \tfrac{\mathbb{E} s^4}{\sigma^2} S \big( X, \Sigma \big)^2 + \tfrac{\mathbb{E} s^4 - \sigma^4}{\sigma^2} \big\lVert \Sigma^{\frac{1}{2}} \big\rVert_{HS}^2
\end{align*}
and we simplify $\big\lVert \Sigma^{\frac{1}{2}} \big\rVert_{HS}^2 = \text{Tr} \ \Sigma = \mathbb{E} \lVert X \rVert^2$.

\hfill $\square$
\end{mdframed}

\begin{corollary}\label{corl_clt}
Suppose $X_i$ are iid with Stein kernel $\tau$ and covariance $\Sigma$, and let $\bar{X} = \tfrac{1}{\sqrt{n}} \sum_{i=1}^n X_i$ be the normalized sum. Then we get a quantitative central limit theorem by combining lemmas \ref{lem_distancefromdiscrepancy}, \ref{lemscalingkernel} (for constant $s = \tfrac{1}{\sqrt{n}}$), and \ref{lemadditivekernel}
\begin{equation*}
\mathcal{W}_2 \big( \bar{X}, \mathcal{N} (0, \Sigma) \big) \leq S \big( \bar{X}, \Sigma \big) \leq \tfrac{1}{\sqrt{n}} S \Big( \textstyle\sum_{i=1}^n X_i, n \Sigma \Big) \leq \tfrac{1}{\sqrt{n}} S \big( X_i, \Sigma \big)
\end{equation*}
\end{corollary}

\subsection{Wasserstein bound in non-isotropic case}\label{apdx_wassersteinbound}

\begin{lemma}\label{lem_distancefromdiscrepancy}
Suppose that $\tau$ is a Stein kernel for random variable $X$, and $\Sigma$ is a symmetric positive-definite matrix. Then
\begin{equation*}
\mathcal{W}_2(X, \mathcal{N}(0, \Sigma)) \leq S(X, \Sigma) \qquad\text{where}\qquad S(X, \Sigma)^2 = \mathbb{E} \lVert \Sigma^{-\frac{1}{2}} \big( \tau(X) - \Sigma \big) \rVert_{HS}^2
\end{equation*}
\end{lemma}

This proof is a compilation of Proposition 3.1 from \cite{Ledoux_HSI} and Lemma 2 from \cite{Otto_Generalization}, additionally keeping track of the covariance matrix. It is based on interpolation of the heat flow along the Ornstein-Uhlenbeck semigroup.

\begin{mdframed}[style=proof]

Let $\mu_0, \mu_\infty$ be measures/densities of $X, \mathcal{N}(0, \Sigma)$ respectively, living in $D$-dimensional space. We will tackle the case when $X$ has a Radon-Nikodym derivative $h = \frac{d \mu_0}{d \mu_\infty}$ with respect to the target normal measure. The general case follows by an approximation argument -- see \cite{Otto_Generalization}.

\mdfsubtitle{HEAT FLOW SEMIGROUP}
Introduce
\begin{equation*}
X_t = e^{-t} X + \sqrt{1 - e^{-2t}} \mathcal{N}(0, \Sigma)
\end{equation*}

Let $\mu_t$ be the measure of $X_t$ and $h_t = \frac{d \mu_t}{d \mu_\infty}$. Define a vector field
\begin{equation}\label{eqdiffusionvf}
v_t(x) = \mathbb{E} \left[ \frac{d X_t}{d t} \Bigg\vert X_t = x \right] \qquad \text{or equivalently} \qquad v_t = - \Sigma \big( \nabla \log h_t \big)
\end{equation}

Then the density $\mu_t$ satisfies the diffusion equation
\begin{equation}\label{eqdiffusiondensity}
\frac{\partial \mu_t}{\partial t} = - \nabla \dotp ( \mu_t v_t )
\end{equation}

A brute-force way to verify the equivalence of definitions in \ref{eqdiffusionvf} and the diffusion equation \ref{eqdiffusiondensity} is to plug in the explicit formulas
\begin{align}
\mu_{0, t}(x_0, x_t) =& (1 - e^{-2t})^{-\frac{D}{2}} \mu_0(x_0) \mu_\infty \left( \frac{x_t - e^{-t} x_0}{\sqrt{1 - e^{-2t}}} \right) \nonumber\\
\mu_t(x) =& (1 - e^{-2t})^{-\frac{D}{2}} \int \mu_0(y) \mu_\infty \left( \frac{x_t - e^{-t} x_0}{\sqrt{1 - e^{-2t}}} \right) dy \label{eqmutexplicit}\\
v_t(x) =& \mu_t(x)^{-1} (1 - e^{-2t})^{-\frac{D}{2}} \int \mu_0(y) \mu_\infty \left( \frac{x_t - e^{-t} x_0}{\sqrt{1 - e^{-2t}}} \right) \left( \frac{-e^{-2t} x + e^{-t} y}{1 - e^{-2t}} \right) dy \nonumber
\end{align}

As a consequence of the diffusion equation \ref{eqdiffusiondensity}, the density $\mu_t$ is transported along the trajectories tangent to $v_t$. Intuitively, if the norm of $v_t$ is small, then the density needs to ``travel a short distance'' to move from $\mu_0$ to $\mu_\infty$. Formally, lemma 2 from \cite{Otto_Generalization} states
\begin{align}
\frac{d^+}{d s} W_2( \mu_t, \mu_{t + s}) \leq& \sqrt{\mathbb{E}  \lVert v_t (X_t) \rVert^2} = \sqrt{\int \lVert v_t(x) \rVert^2 d \mu_t(x)} \nonumber\\
W_2( \mu_0, \mu_\infty ) \leq& \int_0^\infty \sqrt{\mathbb{E}  \lVert v_t (X_t) \rVert^2} dt \label{eq_distfromflownorm}
\end{align}
In the next part of the proof we bound the flow norm as
\begin{equation}\label{eq_flownormbound}
\sqrt{\mathbb{E} \lVert v_t(X_t) \rVert^2} \leq \frac{e^{-2t}}{\sqrt{1 - e^{-2t}}} S(X, \Sigma)
\end{equation}
Substituting this to the inequality~\ref{eq_distfromflownorm} and integrating completes the proof of the lemma.

This was all we need to construct the flow. In the next part of the proof we will need a few more properties. We start with changes in expectations under the semigroup. Define
\begin{equation*}
P_t f (x) \defeq \mathbb{E} f \big( e^{-t}x + \sqrt{1 - e^{-2t}} \mathcal{N}(0, \Sigma) \big) = \int f \big( e^{-t} x + \sqrt{1 - e^{-2t}} y \big) d \mu_\infty(y)
\end{equation*}
This is called Mehler's formula. It is straightforward to check $P_s P_t = P_{s + t}$ and $\mathbb{E} f(X_t) = \mathbb{E} P_t f(X_0)$. Gaussian integration by parts gives a PDE $\frac{\partial}{\partial t} P_t f = \mathcal{L} P_t f$, where $\mathcal{L} = \sum_{i, j} \Sigma_{i j} \partial_i \partial_j - x \dotp \nabla$. Combining these yields
\begin{equation*}
\int f \dot{\mu}_t dx = \frac{d}{ds} \mathbb{E} P_s f(X_t) \Bigg\vert_{s = 0} = \mathbb{E} (\mathcal{L} f) (X_t) = \int f \left( \sum_{i, j} \Sigma_{i j} \partial_i \partial_j + x \dotp \nabla + D \right) \mu_t dx
\end{equation*}

Thus we must have $\frac{\partial \mu_t}{\partial t} = \left( \sum_{i, j} \partial_i \partial_j + x \dotp \nabla + D \right) \mu_t$, which turns out to be a restatement of \ref{eqdiffusiondensity}.

An explicit calculation of $P_t h_0$ turns out to be equivalent to $\frac{\mu_t}{\mu_\infty}$ from formula \ref{eqmutexplicit}, so $h_t = P_t h_0$. It also satisfies
\begin{equation*}
\int f \cdot P_t g d \mu_\infty = \mathbb{E} \left[ f(x) g(y) \Bigg\vert \begin{pmatrix} X \\ Y \end{pmatrix} \sim \mathcal{N} \left( 0, \begin{pmatrix} \Sigma & e^{-t} \Sigma \\ e^{-t} \Sigma & \Sigma \end{pmatrix} \right) \right] = \int P_t f \cdot g d \mu_\infty
\end{equation*}
and $\nabla P_t f = e^{-t} P_t \nabla f$. The diffusion operator satisfies $\int f \mathcal{L} g d \mu_\infty = - \int (\nabla f)^\top \Sigma (\nabla g) d \mu_\infty$.

\mdfsubtitle{BOUND ON THE FLOW NORM}

This part of the proof is concerned with proving the inequality~\ref{eq_flownormbound}. We start from the transformations
\begin{align*}
\int \lVert v_t \rVert^2 d \mu_t =& \int (\nabla \log h_t)^\top \Sigma^2 (\nabla h_t) d \mu_\infty = \\
=& e^{-t} \int (\nabla \log h_t)^\top \Sigma^2 (P_t \nabla h_0) d \mu_\infty = \\
=& \int (\nabla P_t \log h_t)^\top \Sigma^2 (\nabla h_0) \mu_\infty(x) dx = \\
=& - \int \nabla \dotp \left( \mu_\infty \cdot \Sigma^2 \nabla P_t \log h_t \right) h_0 dx = \\
=& \int \left( x \dotp \Sigma \nabla P_t \log h_t - \nabla \dotp \Sigma^2 \nabla P_t \log h_t \right) h_0 \mu_\infty dx = \\
=& \int \left( x_i \Sigma_{i j} \partial_j P_t \log h_t - \Sigma_{i k} \Sigma_{k j} \partial_i \partial_j P_t \log h_t \right) d \mu_0 = \\
=& \int \big( \tau_{i k}(x) - \Sigma_{i k} \big) \Sigma_{k j} \partial_i \partial_j P_t \log h_t d \mu_0(x)
\end{align*}

where in the last two lines we used the Einstein summation convention. We substitute the identity
\begin{equation*}
\partial_i \partial_j P_t \log h_t = \frac{e^{-2t}}{\sqrt{1 - e^{-2t}}} \int \big( \Sigma^{-1} y \big)_i \big( \partial_j \log h_t \big) \left( e^{-t} x + \sqrt{1 - e^{-2t}} y \right) d \mu_\infty(y)
\end{equation*}
to get
\begin{equation}\label{eq_flownorm_raw}
\mathbb{E} \lVert v_t(X_t) \rVert^2 = \frac{e^{-2t}}{\sqrt{1 \! - \! e^{-2t}}} \iint \left( y^\top \Sigma^{-1} \big( \tau(x) \! - \! \Sigma \big) \right) \cdot \big(\Sigma \nabla \log h_t \big) \left( e^{-t} x + \sqrt{1 \! - \! e^{-2t}} y \right) d \mu_0(x) d \mu_\infty(y)
\end{equation}
By Cauchy-Schwarz inequality the integral is at most
\begin{align}
& \sqrt{\iint \left\lVert y^\top \Sigma^{-1} \big( \tau(x) - \Sigma \big) \right\rVert^2 d \mu_0(x) d \mu_\infty(y)} \times \label{eqsteinkernelroot}\\
&\times \sqrt{\iint \left\lVert \big( \Sigma \nabla \log h_t \big) \left( e^{-t} x + \sqrt{1 - e^{-2t}} y \right) \right\rVert^2 d \mu_0(x) d \mu_\infty(y)} \label{eqflownormroot}
\end{align}

Expression under the root in \ref{eqsteinkernelroot} equals
\begin{align*}
\iint y_i y_j & \Sigma^{-1}_{i k} \Sigma^{-1}_{j l} \big( \tau(x) - \Sigma \big)_{k m} \big( \tau(x) - \Sigma \big)_{l m} d \mu_0(x) d \mu_\infty(y) = \\
=& \int \Sigma_{k l}^{-1} \big( \tau(x) - \Sigma \big)_{k m} \big( \tau(x) - \Sigma \big)_{l m} d \mu_0(x) = \\
=& \int \left\lVert \Sigma^{-\frac{1}{2}} \big( \tau(x) - \Sigma \big) \right\rVert_{HS}^2 d \mu_0(x) = S (X, \Sigma)^2
\end{align*}
while the expression under the root in \ref{eqflownormroot} is
\begin{equation*}
\int P_t \left( \lVert v_t \rVert^2 \right) d \mu_0 = \int \lVert v_t \rVert^2 \cdot P_t h_0 d \mu_\infty = \int \lVert v_t \rVert^2 d \mu_t
\end{equation*}

These two simplifications allow to bound the equation~\ref{eq_flownorm_raw} as
\begin{equation*}
\mathbb{E} \lVert v_t (X_t) \rVert^2 \leq \tfrac{e^{-2t}}{\sqrt{1 - e^{-2t}}} \cdot S \big( X, \Sigma \big) \cdot \sqrt{\mathbb{E} \lVert v_t (X_t) \rVert^2}
\end{equation*}
Which is equivalent to the inequality~\ref{eq_flownormbound}. Now, combining inequalities~\ref{eq_distfromflownorm} with~\ref{eq_flownormbound} completes the proof of the lemma.

\hfill $\square$
\end{mdframed}

\end{document}